\documentclass[12pt]{article}
\usepackage{amssymb}
\def \C{{\Bbb C}}

\def \Z{{\Bbb Z}}
\def \p{{\Bbb P}}

\begin{document}
\title {Exceptional points of an endomorphism of the projective plane} 
\author{E. Amerik and F. Campana} 
\date{October 15, 2003}

\maketitle

{\bf 1. Introduction}

\

Let $f:\p^k\rightarrow \p^k$ be an endomorphism of the complex 
projective space, of degree bigger than one. An algebraic subset $V\in \p^k$
is called {\it completely invariant} for $f$, if
$f^{-1}(V)=V$. The {\it exceptional set} of $f$ is the largest completely
invariant proper algebraic subset of $\p^k$ (it is observed for example in
\cite{BD} that this definition does make sense). 

It seems to be interesting from the point of view of
dynamics to understand which algebraic subsets of $\p^k$ can actually be
exceptional or completely invariant. 
Recently, this question was studied by several authors
(\cite{FS}, \cite{CLN}, \cite{BCS}). For example, in \cite{BCS} it is
shown that an exceptional set is always a union of linear subspaces.
It is easy to see
(and observed in \cite{FS}) that there are at most $k+1$ {\it hyperplanes}
in an exceptional set. This bound is of course sharp: take $f$ raising
the homogeneous coordinates to $m$th power. 
However, no such bound was found, up to now,
for the number of higher-codimensional components of an exceptional set.
In this article, we first treat the simplest case of this problem: that
of points in $\p^2$, and then show how to obtain a bound for the number
of completely invariant linear subspaces of codimension two.

More precisely, let $f:\p^2\rightarrow \p^2$ be an endomorphism of degree
bigger than one. Let us say that $f$ is {\it completely ramified} at 
a point $P\in \p^2$, if $f^{-1}(P)$ consists of one point. In Section 3, 
we prove the
following

\

{\bf Theorem 1:} 
{\it $f$ can be completely ramified at nine points at most.
No three such points are on a line and no six on a conic.
Moreover, if $f$ is completely ramified at nine points, 
they lie on a unique cubic, which
has
a singularity at one of them.}

\

Our result does not seem to be the best possible. Indeed, we do not know
any example of an $f$ which is completely ramified at more than three points
(an obvious $f$ completely ramified at three points is given by the formula
$$f(x:y:z)=(x^m:y^m:z^m);$$ one can also take maps of the type
$$f(x:y:z)=(x^m:y^m:z^m+xyg(x,y,z))$$ and compositions of such).

Let us now explain the
idea of the proof of this theorem.
 Let $m$ be such that $f^*({\cal O}_{\p^2}(1))=
{\cal O}_{\p^2}(m)$, so that $deg(f)=m^2$. Suppose that $f$ is completely
ramified at $3d$ points $P_1, \dots , P_{3d}$ 
lying on a smooth curve $D$ of degree $d$, and
consider the (rarely occuring of course) case when the inverse image $C$ of
$D$ is smooth and reduced. By adjunction, we have $deg(K_C)=md(md-3)$.
On the other hand, at each $P_i$, the ramification index of 
$f|_C:C\rightarrow D$ is $m^2-1$, and by Hurwitz' formula, $$deg(K_C)
\geq deg((f|_C)^*K_D)+3d(m^2-1),$$
so that 
$$md(md-3)\geq m^2d(d-3)+3d(m^2-1),$$
and this is impossible if $m>1$.

Now through any nine points in the plane there is a cubic. This cubic, and
a fortiori its inverse image, does not have to be smooth. But it happens
that one can, with some modifications, make the above argument work
for singular $D$ and $C$, replacing the canonical class by the dualizing
sheaf. One obtains then that there are at most $3d-1$ points of complete
ramification in the {\it smooth} locus of any curve of degree $d$ in $\p^2$,
and also at most $3d-1$ points of complete ramification on any curve
which does not contain a component of the branch locus. This implies almost
immediately that $f$ can be completely ramified at 11 points at most, and 
some further analysis
implies our theorem. An important point of this analysis is that a
point
of complete ramification must be of very high multiplicity on the
direct image of the ramification divisor (see Proposition 1 below).

In the rest of the paper, we generalize this to higher dimension. 
We obtain

\

{\bf Theorem 2:} 
{\it Let $f:\p^{N-1} \rightarrow \p^{N-1}$ be a morphism of degree
bigger than one. The number of  
codimension-two completely invariant subspaces for $f$ is less than $4N^2$.}

\

The proof of this theorem goes essentially along the same lines. One considers
a general $\p^2$ in $\p^{N-1}$, the intersection points $P_i$ 
of our completely
invariant subspaces with this $\p^2$ and the restriction of $f$ 
to the inverse image
of this $\p^2$. One obtains this time (Proposition 3) that a curve
of degree $d$ contains at most $Nd-1$ of the points $P_i$ in its smooth locus,
and a ``flexible'' curve of degree $d$ contains at most $Nd-1$ of the
points $P_i$. Again, this implies almost immediately a somewhat weaker
bound $6N^2$. We can improve this to $4N^2$, either by recalling  
a Cayley-Bacharach type statement due to Ellia and Peskine (a remark in the end
of the section 4), or by proving a general combinatorial result
 on
points of high multiplicity on an effective divisor in the plane (Section 5),
which in particular implies our Theorem 2. Again, our bound is certainly
not the best possible; it is however clear that $N^2$ must appear, 
as the map raising the
coordinates to the $m$th power has $\frac{1}{2}N(N-1)$ completely invariant
subspaces of codimension 2. 

We are grateful to N. Sibony, S. Cantat,
T.-C. Dinh and C. Favre for asking the question, and to A. Otwinowska and
D. Perrin for providing us with some useful references on the 
Hilbert function of a set of
points in the plane.

\

{\bf 2. Some remarks on differential forms}

\

Let us first recall a few basic facts on the dualizing sheaf of a curve
on a smooth surface (see for example \cite{BPV}, chapter II), 
and make some observations on the behaviour of this 
sheaf under morphisms.

Let $D\subset S$ be an effective divisor on a smooth surface $S$. The
dualizing sheaf $\omega_D$ is then the line bundle ${\cal O}(K_S+D)|_D$.
If $D$ is reduced, then $\omega_D$ identifies with a sheaf of
meromorphic differentials on the normalization $\tilde{D}$, called
Rosenlicht differentials. Concretely, these are the differentials
$\phi$ such that for any $x\in D$ and $g\in {\cal O}_{D,x}$, one has
$$\sum_{y\in n^{-1}(x)}res(y,g\phi)=0,$$ where $n$ is the normalization
map.
Moreover one can describe $\omega_D$ in local
coordinates in the following way: if a local equation of $D$ in $S$
is $F(x,y)=0$, and $x$ is non-constant on any local branch of $D$,
then $\omega_D$ is locally generated by $\frac{dx}{F'_y}$ (or, more
precisely, by $n^*(\frac{dx}{F'_y})$),
which one views as the ``residue'' of $\frac{dx\wedge dy}{F}$.

\

{\bf Lemma 1:} {\it Let $f: X\rightarrow Y$ be a finite morphism of smooth 
surfaces, $D\subset Y$ a (reduced) curve and $C=f^*D$ the (scheme-theoretic)
inverse image of $D$. Then $f$ induces the injection of sheaves
$(f|_{C_{red}})^*\omega_D
\rightarrow \omega_{C_{red}}$.}

\

{\it Proof:} Let $C=\sum_{i=1}^l n_iC_i$ be the decomposition of the divisor $C$ into
the sum of its reduced irreducible components, 
so that the ramification divisor of
$f$ is $R=\sum_{i=1}^l (n_i-1)C_i + R'$, where $R'$ contains no component of $C$.
The map $f$ induces the injection  $f^*K_Y\rightarrow K_X$, locally given
by the multiplication by the Jacobian. The image of $f^*K_Y$ is of course
contained in the subsheaf of forms vanishing on $C_i$ with multiplicity
$n_i-1$, so we get an injection 
$ f^*K_Y\rightarrow K_X-\sum_{i=1}^l (n_i-1)C_i$. Adding $C$
and restricting to  $C_{red}$, we get a morphism $(f|_{C_{red}})^*\omega_D
\rightarrow \omega_{C_{red}}$, which is still an injection as $C$ and $R'$
do not have common components. 

\

{\bf Remarks:} 1) It is clear that in general, a morphism $h:C\rightarrow D$,
even a surjective one,
does not induce a morphism of the dualizing sheaves when $D$ is not
smooth (take for example
the normalization map of a cuspidal curve).

2) The cokernel of our injection is the structure sheaf of the intersection
of $R'$ and $C$, so that the contribution of each ``ramification point''
to the difference of $deg(\omega_{C_{red}})$ and $deg((f|_{C_{red}})^*\omega_D)$
is the intersection index of $R'$ and $C$ at that point.

3) One could also write all this in local coordinates
using the above explicit description of the dualizing sheaf of a reduced 
curve.

4) Obviously we also have a morphism from $(f|_C)^*\omega_D$ to $\omega_C$,
and it is not injective when $C$ is non-reduced (neither is its cokernel
concentrated at points). 

\

We shall need a local computation in order to estimate the ``ramification
multiplicities'' at the points of complete ramification. Let now
$(D,0)$ be a germ of a smooth curve, $z$ a coordinate on $D$. Let
$(C,(0,0))$ be a germ of a reduced curve on a smooth surface, given 
by the equation $F(x,y)=0$. Let further $n:\tilde{C}\rightarrow C$ be the normalization
of $C$, and $Q_1, \dots , Q_r$ the points of $\tilde{C}$ with $n(Q_i)=(0,0)$.

\

{\bf Lemma 2:} {\it Let $g:(C,(0,0))\rightarrow (D,0)$ 
be a germ of a finite map
of local degree $\delta$ at $(0,0)$. Denote by $\mu_i$ the
vanishing order at $Q_i$ of $n^*g^*dz$ as a section of $n^*\omega_C$.
Then $\sum_{i=1}^r \mu_i \geq \delta -1$.}

\

{\it Proof:} 
Let $e$ be the degree of the restriction to $C$ of the first
projection $p_1:(x,y)\rightarrow x$. We may suppose that the function
$F(x,y)$ defining $C$ is in the Weierstrass form:

$$F(x,y)=y^{e}+a_{e-1}(x)y^{e-1}+\dots +a_1(x)y+a_0(x),$$

with $a_i(0)=0$ for any $i$.

Suppose first that $C$ is an irreducible germ. Choose a
coordinate $t$ on $\tilde{C}$ so that the map $p_1\cdot n$ is given
by $x=t^e$ and the map $g\cdot n$ is given by $z=h(t)t^{\delta}$, where 
$h$ is holomorphic with $h(0)\neq 0$. We have then $n^*dx=et^{e-1}dt$
and $n^*g^*dz=h_1(t)t^{\delta -1}dt$ with $h_1$ holomorphic.

As $n^*\omega_C$ is generated by $n^*(\frac{dx}{F'_y})$, all that we have 
to show is that $n^*g^*dz=h_2(t)t^{\delta -1}n^*(\frac{dx}{F'_y})$ for some
holomorphic $h_2$, i.e. that $n^*F'_y$ has a zero of order at least
$e-1$ at $0$. But this is obvious, as $n^*y$ has a zero at $0$, and
$n^*x$ has a zero of order $e$.

Now if $C$ is reducible, let $C_j, j=1, \dots, l$ be the irreducible
components and  $F(x,y)=\Pi_{j=1}^l F_j(x,y)$, each $F_j$ defining $C_j$.
Denote by $n_j:\tilde{C_j}\rightarrow C_j$ the normalization and by $\delta_j$
the degree of $g|_{C_j}$. As $\delta$ is the sum of $\delta_j$, it suffices to
show that $g^*n_j^*dz$ has a zero of order at least $\delta_j$ at $n_j^{-1}(0,0)$,
again as a section of $n^*\omega_C$, of course. 

This is the same
kind of computation as before:
choose appropriate coordinates, write
$$F_j(x,y)=y^{e_j}+a_{j,e_j-1}(x)y^{e_j-1}+\dots +a_{j,1}(x)y+a_{j,0}(x),$$
$e_j$ being the degree of $p_1\cdot n_j$. Repeating the previous 
considerations, we see that
this time we have to show that $n_j^*F'_y$ has a zero of order at
least $e_j$ at $0$. But $F=F_jG$ with $G$ holomorphic and vanishing
at $(0,0)$ and $F'_y=(F_j)'_yG+F_jG'_y$. We have that $n_j^*(F_j)'_y$ has
a zero of order at least $e_j-1$, $n_j^*F_j$ is identically zero,
 and $n_j^*G$ has a zero at $0$, which proves our assertion.

\

This lemma, together with the projection formula, implies the following
fact which will be quite useful:

\

{\bf Proposition 1:} {\it  Let $f:X\rightarrow Y$ be a finite (proper) 
morphism between
smooth surfaces (or their germs). Let $R$ denote its ramification divisor.
For any $y\in Y$, the following inequality holds:}

$$mult_y(f_*R)\geq \sum_{x\in f^{-1}(y)}(\delta_x(f)-1),$$
{\it where $mult_y$ stands for the multiplicity at $y$ and $\delta_x(f)$
denotes the local degree of $f$ at $x$.}

\

{\it Proof:} Choose a sufficiently general smooth curve $D$ through $y$, 
so that $mult_y(f_*R)$ equals to the local intersection index 
$(f_*R\cdot D)_y$. By the projection formula, this intersection index is
equal to $\sum_{x\in f^{-1}(y)}(R\cdot f^*D)_x$. It remains to apply
Lemma 2 together with the second remark after Lemma 1.

\

Finally, we prove one more lemma about differentials on singular curves,
though it is not absolutely necessary for the sequel:

\

{\bf Lemma 3:} {\it Let $g:(\C^2,0)\rightarrow (\C^2,0)$ be a germ of a
finite holomorphic map, of local degree $\delta$. Let $(D,0)\subset (\C^2,0)$
be a germ of a (reduced) curve, such that $(C,0)=g^*(D,0)$ is reduced, and let
$r$ be the multiplicity of D at 0. Then the length of the cokernel of the lift
to the normalization of $C$ of
the injection of the dualizing sheaves $(g|_C)^*\omega_D\rightarrow \omega_C$ is
at least $r(\delta-1)$.}

\

{\it Proof:} Since $C$ is reduced, this length is the intersection
multiplicity of $C$ and the ramification divisor $R$, 
or, by projection formula,
the intersection multiplicity of $D$ and $f_*R$. This multiplicity
is at least equal to the product of the multiplicities at $0$ of
$D$ and $f_*R$, so the lemma follows from Proposition 1. 
 
\

{\bf Remark:} The condition that $D$ is smooth is essential for Lemma
2 ( moreover, Lemma 2 does not make sense for a singular $D$, since
in this case the inverse image of a Rosenlicht differential on $D$ does not
have to be a Rosenlicht differential on $C$ ). 
Lemma 3 as it is also does not make sense for non-reduced
$C$, for example, because in this case the cokernel is not concentrated 
in points, and also because while a general point of $C$ should not contribute
at all to the difference of the degrees of $(g|_C)^*\omega_D$ and $\omega_C$,
the local degree  $\delta$ at such a point is still bigger than one
(this however should be possible to repair by considering the local degree
of $g|_{C_{red}}$ rather than $\delta$).
We do not know if there is some analogue of Lemma 3 for non-reduced curves, and
it is an interesting question (a suitable analogue could simplify and
improve the proof of our theorems).
In any case there are many examples showing that if $D$ is singular
and
$C$ is non-reduced, the length of the cokernel of
$(g|_{C_{red}})^*\omega_D\rightarrow \omega_{C_{red}}$ need not be bounded from below
by a non-trivial function of the local degree of the map from
$C_{red}$ to $D$. We list below several of them, in which
this length is zero, because the ``residual ramification'' $R'$ is 
empty:

\

{\bf Examples:}

1) $g$ given by the formula $u=x^m,v=y^m$, and $D=\{(u,v):uv=0\}$.
The ramification divisor $R$ is given by the equation $x^{m-1}y^{m-1}=0$,
and $C_{red}$ by the equation $xy=0$. Our inclusion 
$(g|_{C_{red}})^*\omega_D\rightarrow \omega_{C_{red}}$ is thus an isomorphism.

2) An example with locally irreducible $D$: let $m>1$ be an odd
integer and let
$g$ be given by
$u=(x^m-y^m)/2, v=xy$. 
The equation of the ramification divisor
is then $x^m+y^m=0$. Consider $D=\{(u,v): u^2+v^m=0$. The equation
of the inverse image of $D$ is $((x^m-y^m)/2)^2+(xy)^m=(x^m+y^m)^2/4=0$,
so $g^*(D)=2R$ and $\omega_{C_{red}}=(g|_{C_{red}})^*\omega_D$.

3) A more exotic example with cuspidal $D$ is the following: let $V\subset 
{\Bbb C}^3$ be the canonical hypersurface singularity given by the equation 
$u^2+v^3+w^5=0$.
It is a quotient singularity $h:{\Bbb C}^2\rightarrow V=({\Bbb 
C}^2/A_5)$, for a certain $Sl_2({\Bbb C})$-representation of $A_5$. 
The map $h$ is unramified except at the origin.
Let $p(u,v,w):=(u,v)$, and $g:{\Bbb C}^2\rightarrow {\Bbb C}^2$ be 
the composition $p\circ h$. Let $D$ be given by the equation 
$u^2+v^3=0$. Let $C=g^*(D)$,
and let $R$ be the ramification divisor of $g$. One checks that 
$R=4C_{red}$ and 
$g^*(D)=5C_{red}$, that is, $\omega_{C_{red}}=(g|_{C_{red}})^*\omega_D$.
 In this example, the degree of $g$ is 300, and the
degree
of $g$ restricted to $C_{red}$ is 60. We do not know if there are
examples with cuspidal $D$ and arbitrarily high degree of $g|_{C_{red}}$.

\

{\bf 3. Proof of Theorem 1}

\

We start with the following

\

{\bf Proposition 2:} {\it Let $f$ be an endomorphism of projective plane,
of degree bigger than one. Let $D\subset \p^2$ be a curve of degree $d$.
Then $f$ cannot be completely ramified at $3d$ (or more) 
points of the smooth locus
of $D$.}

\

{\it Proof:} Let $f$ be an endomorphism as above, completely ramified at $3d$ 
smooth points of $D$. Let $f^*{\cal O}(1)={\cal O}(m)$, let  $C=f^*D$ 
and let $c$
be the degree of $C_{red}$. By the projection formula, the degree of
the restriction of $f$ to $C_{red}$ is $\frac{mc}{d}$. We have, by Lemma 1,
the induced injection 
$\alpha: (f|_{C_{red}})^*\omega_D\rightarrow \omega_{C_{red}}$, and,
as we suppose that our $3d$ points of complete ramification are
in the smooth locus of $D$, by Lemma 2, each of them contributes at least $\frac{mc}{d}-1$
to the length of the cokernel of the lift of this injection 
to the normalization of $C_{red}$. So that we have:

$$deg(\omega_{C_{red}})\geq 3d(\frac{mc}{d}-1) + \frac{mc}{d}deg(\omega_D),$$

that is,

$$c(c-3)\geq dmc-3d,$$ which of course implies that $c\leq d$, since by the 
definition of $c$ we have $dm\geq c$, and even that $c<d$, since $m>1$. 
If moreover $mc-3>0$ we conclude
that $c(c-3)\geq c(mc-3)>c(c-3)$, a contradiction. So that we must have
$m=2$ and $c=1$. Then $d=2$, since $\frac{mc}{d}$ is an integer
(the degree of the restriction of $f$ to $C_{red}$). 
This implies that the restriction of $f$ to $C_{red}$ is one-to-one
and so $C=4C_{red}$. Since the ramification divisor of $f$ is of degree three,
it must be $3C_{red}$, so that the complement to $C_{red}$ in $\p^2$ is
an unramified covering of degree four of the complement to $D$ in $\p^2$.
This is a contradiction because the fundamental group of the latter is
$\Z/2\Z$.

\

Notice that this proposition already implies   
that an endomorphism $f$ (of degree $>1$) of projective plane cannot be 
completely
ramified at more than 11 points, by the following general remark:

\

{\bf Lemma 4:} {\it Let $E$ be a subset of the projective plane, such that
any curve of degree $1\leq d\leq 4$ contains at most $3d-1$ points of $E$
in its smooth locus. Then $E$ cannot contain more than 11 points.}

\

{\it Proof:} Suppose that we have 12 points 
$P_1, \dots, P_{12}$ in $E$. It is immediate from the assumption
that no three of those points are on a line, no six on a conic, 
no ten on a cubic,
and through any nine there is a unique cubic which must be 
singular at one of the $P_i$
(but irreducible). Consider the linear system of quartics passing through
$P_1, \dots, P_{12}$. It is a linear subspace $L$ of dimension at least two
in the projective space $P=\p^{14}$ parametrizing all the plane quartics.
Suppose that all these quartics are singular. Then, by Bertini, they are
all singular at one of the $P_i$, say $P_1$. So $L$ is contained in the
linear subspace $M\subset P$ parametrizing the quartics singular at $P_1$.
The space $M$ is of dimension 11. Now consider the subset $V$ of $M$ 
parametrising reducible quartics which are of the form
$Q=H\cup C$ with $H$ a line through $P_1$ and $C$ is a (possibly reducible)
cubic through $P_1$. Clearly $V$ is a subvariety of $M$ of dimension 9.
So $L$ and $V$ must intersect; that is, among our quartics passing through
12 points, there is one which is a union of a line and a cubic. As
we have already seen that no three points are on a line and no ten on a cubic,
this gives a contradiction. So among our quartics there must be a
smooth one, but this again contradicts our assumption.

\

The proof of Theorem 1 does not depend on Lemma 4:

\

{\it Proof of Theorem 1:} Consider first an endomorphism $f$ with
nine points of complete ramification $P_1, \dots, P_9$. We see from
 Proposition 2 that through these points there is a unique irreducible
cubic $D$ which has a singularity
at one of our points, say, $P_1$. Let $R$ denote, as usually, the
ramification divisor. It is of degree $3m-3$, and its direct image $f_*R$
is, by the projection formula, of degree $3m^2-3m$. Write $f_*R=aD+D'$
where $a\geq 0$ and $D'$ does not contain $D$ as a component. 
We shall show
that the number $a$ is rather big.

Indeed, let us compute
the intersection number
$D'\cdot D$: on the one hand, it is equal to
$$(f_*R-aD)\cdot D=9m^2-9m-9a,$$
on the other hand, by Proposition 1, the multiplicity of $D'$ at $P_1$
is at least $m^2-1-2a$ and the multiplicity of $D'$ at $P_i$, $2\leq i\leq 9$,
is at least $m^2-1-a$. The intersection number is then at least
$8(m^2-1-a)+2(m^2-1-2a),$ so that
$9m^2-9m-9a\geq 10m^2-10-12a,$
or, $$a\geq\frac{m^2+9m-10}{3}.$$ 

Now if we have ten points in which $f$ is completely ramified, they are
not all on a cubic, so there are ten cubics $D_1, \dots, D_{10}$, each $D_i$
containing exactly nine points of our ten, and each being a component
of $f_*R$ of multiplicity at least $\frac{m^2+9m-10}{3}$.
But since $10(m^2+9m-10)>3m^2-3m$, this is impossible.

\

{\bf Remark:} We see that we have not completely generalized the argument
from the Introduction: that is, we have not found a ``Hurwitz formula''
(that is, a statement similar to Lemma 2 and Lemma 3)
for differentials on a non-reduced curve and so in our proof
we have to circumvent the case where $C=f^*D$ is non-reduced. If $D$ is smooth,
we do it by considering $C_{red}$, but if it is not, then there are examples
where the contribution of a point of complete ramification to 
$deg(\omega_{C_{red}})$ is not sufficient for our purposes (for
instance,
zero, as in our examples from the previous section). Anyway, even if 
we could understand those differentials better, this would only give an upper
bound of 8 points. It is of course impossible to prove by our method
that there are at most three. On the other hand, this could considerably
improve our Theorem 2, giving immediately a bound $2N^2$, whereas
in what follows we obtain rapidly $6N^2$ and can lower it to $4N^2$ only
with some effort.

\

\

{\bf 4. Generalizing to higher dimension}

\

Let
$f:\p^{N-1}\rightarrow \p^{N-1}$, $N>3$,
 be a morphism of degree bigger than one,
and $f^*{\cal O}(1)={\cal O}(m)$. Let $L_1, \dots L_k$ be codimension-two
subspaces of $\p^{N-1}$ such that for every $1\leq i\leq k$, 
$f^{-1}(L_i)= L_i$. Notice that, contrary to the situation of the last
section, we suppose that our linear subspaces are completely
invariant (rather than have another linear subspace as the inverse image).
This is because we will need to assume that $m$ is sufficiently
big with respect to $N$ ( compare with the proof of Proposition 2 ).
If we want to estimate the number of completely invariant linear
subspaces, we may of course assume this, iterating $f$ if necessary.

 Consider a general $S=\p^2\subset \p^{N-1}$.
By a Bertini-type theorem, the inverse image $T$ of this $\p^2$ is a
smooth surface; obviously it is a complete intersection of type 
$(m,m,\dots ,m)$, of degree $m^{N-3}$. Let $P_i$ be the intersection point of $S$ and $L_i$. The degree of the map from $f^{-1}(L_i)$ to $L_i$ is 
$m^{N-3}$, and the local degree of $f$ at a general point of $f^{-1}(L_i)$
is $m^2$. We see thus that each $P_i$ has $m^{N-3}$ points 
$Q_{i,j}\in T, j=1, \dots ,
m^{N-3}$, in its inverse image, and the local degree of $f_T$ (the 
restriction
of $f$ to $T$) at those points is $m^2$. In what follows we shall get
a bound for $k$ by studying the map $f_T: T\rightarrow S= \p^2$.

Let us introduce some more notations: 
denote by $E$ the set of $P_i$, $1\leq i\leq k$, and by $H_T$ the hyperplane
section class on $T$. The canonical class $K_T$ of $T$ is thus 
$(-N+(N-3)m)H_T$, and the ramification divisor of $f_T$ is $(m-1)NH_T$
(so that $f_*R$ is of degree $m^{N-2}(m-1)N$).

Here is an analogue of Proposition 2, adapted to our purposes:

\

{\bf Proposition 3:} {\it In the situation as above, let 
$D$ be a (reduced) curve of degree $d$ on $S=\p^2$. If $m>>N$, then

a) $D$ has at most $Nd-1$ points of $E$ in its smooth locus;

b) If $D$ does not contain any component of the branch locus, then
there are at most $Nd-1$ points of $E$ on $D$;

c) There are at most $Nd+\frac{(d-1)(d-2)}{2}-1$ points of $E$ on $D$.}

\

{\it Proof:} a) Suppose that there are $Nd$ or more points of $E$
in the smooth locus of $D$. Let $C$ be 
the (scheme-theoretic) inverse image of $D$
under $f_T$. $C$ is a curve cut out
on $T$ by an hypersurface of degree $md$, and of course it need not be 
reduced; let $C_{red}$ denote its reduction. Its dualizing sheaf 
$\omega_{C_{red}}$ is a line bundle because $C_{red}$ is a divisor
on a smooth surface $T$. We cannot directly
compute the degree of $\omega_{C_{red}}$, because, contrary
to the case $N=3$, we cannot be sure that $C_{red}$ is a complete
intersection. But we can get an inequality from the Hodge index theorem.
Indeed, let
$c= \frac{C_{red}H_T}{m^{N-3}}$ (if $C_{red}$ was a complete intersection, it
would of course mean that
$C_{red}$ is cut out on $T$ by an hypersurface of degree $c$).
By the projection formula, we have that 
the degree of the map from $C_{red}$ to $D$ is $\frac{mc}{d}m^{N-3}$.
Now by the Hodge index theorem, $C_{red}^2H_T^2\leq (C_{red}H_T)^2$,
that is, $C_{red}^2\leq c^2m^{N-3}$.
Finally, we have:
$$deg(\omega_{C_{red}})=K_TC_{red}+C_{red}^2 \leq (-N+c+(N-3)m)cm^{N-3},$$
and of course $$deg(\omega_D)=(d-3)d.$$

Denote now by $e_{i,j}$ the local degree of the map from $C_{red}$ to $D$
at those $Q_{i,j}$ which are on $C_{red}$ 
(if $C$ was reduced, we would of course have $e_{i,j}=m^2$).
Obviously we have:
$$\sum_{j=1}^{m^{N-3}}e_{i,j}=\frac{mc}{d}m^{N-3}$$
for all $i$ such that $P_i$ is on $D$. By Lemma 2, if $P_i$ is a smooth
point of $D$, then
at each of the corresponding points $Q_{i,j}$, 
the length of the cokernel of the injection of the dualizing sheaves is 
at least $e_{i,j}-1$.

This, and the assumption that there are at least $Nd$ points of $E$
in the smooth locus of $D$, implies the following inequality:

$$(\ast)(-N+c+(N-3)m)cm^{N-3}\geq \frac{mc}{d}m^{N-3}d(d-3)+
Ndm^{N-3}(\frac{mc}{d}-1),$$
which reduces to 
$$c^2-Nc\geq mcd-Nd.$$
This inequality, as in the proof of Theorem 1, cannot hold if $mc-N>0$.

As was already said, we assume 
 that the subspaces $L_i$ are completely invariant,
that is, $f^{-1}(L_i)$ is not an arbitrary linear subspace, but $L_i$
itself, and so we may assume that $m$ is arbitrarily big with respect to $N$,
replacing $f$ by a suitable power
of $f$. The number $c$ however is not necessarily an integer and depends
on the map $f$ and the curve $D$, so we need to make an estimation for
$c$ to conclude. In fact it is not difficult to bound $c$ from below 
as follows: the ramification divisor of our original map
$f:\p^{N-1}\rightarrow \p^{N-1}$ is a (possibly reducible and non-reduced) 
hypersurface
of degree $mN-N$. Therefore, at a general point of any component of
this divisor, the local degree of $f$ is at most $mN-N+1$, and the image
of the set of points where this local degree is bigger, is a subvariety
of codimension at least two in $\p^{N-1}$. Since our
surface $S$ is a general plane in $\p^{N-1}$, we can choose it so that
it intersects this subvariety in a finite set of points. For such a 
choice of $S$, and any reduced curve $D\subset S$, the multiplicity of 
any irreducible
component of $C=f^*D$ is at most $mN-N+1$. This gives
$$c\geq \frac{md}{mN-N+1},$$
and, taking for example $m\geq N^2$, we obtain $mc>N$. This proves the part
a).

b) The condition that $D$ does not contain a component of the branch locus
means that $f^*D$ is reduced. So one could repeat a part of the 
above argument applying 
Lemma 3 instead of Lemma 2 and remarking that now $c=md$, so things become
easier. But it is probably better
to apply Proposition 1 directly. Indeed, $f_*R$ is, by projection formula,
of degree $m^{N-2}(m-1)N$. Suppose that $D$ passes through $r$ points of 
complete
ramification. By Proposition 1, these are points of multiplicity at least 
$(m^2-1)m^{N-3}$ on $f_*R$.
One has then $m^{N-2}(m-1)Nd=f_*R\cdot D\geq r(m^2-1)m^{N-3}$, so $r<Nd$.

c) If $D$ is irreducible, this is just a reformulation of a), 
since $\frac{(d-1)(d-2)}{2}$ is the arithmetic genus $p_a(D)$ of $D$, so $D$
cannot have more than $\frac{(d-1)(d-2)}{2}$ singularities.
Otherwise, $D=\bigcup D_i$ with $D_i$ irreducible of degree
$d_i$. Each $D_i$ contains less than $Nd_i+p_a(D_i)$ points of complete
ramification. Their union $D$ contains then less than $Nd+p_a(D)$ 
of those points, as follows from the inequality
 $p_a(D)\geq
\sum p_a(D_i)$ (implied by standard exact sequences).

\

Notice that, as in the previous section, Proposition 3 implies almost
immediately that $|E|<6N^2$:

\

{\bf Proposition 4:} {\it Let $N$ be a natural number and $E$ 
a set of points in the plane, with the following properties:

(1) Any curve of degree $d$ contains at most $Nd-1$ points of $E$ in its
smooth locus;

(2) Any ``flexible'' curve (that is, a general member of
a linear system without base components) of degree $d$
contains at most $Nd-1$ points of $E$.  

Then $E$ has less than $6N^2$ points.}

\

{\it Proof:} Suppose this is not true. 
Fix a subset $V\subset E$ of $6N^2$ points
and consider the linear system of curves of degree $4N$ containing $V$.
The dimension of this linear system is at least $4N(4N+3)/2-6N^2=2N^2+6N$.
Denote by $l$ the degree of the union of its base components, i.e. of its
fixed part. Remark that $l$ cannot be very big: indeed, 
the degree of the free part
is $4N-l$ and so the dimension of our linear system is at most
$(4N-l)(4N-l+3)/2$, so we get an inequality
$(4N-l)(4N-l+3)\geq 4N^2+12N$, or $12N^2-8Nl+l^2\geq 0$, from where
$l\leq 2N$. Now the property (1) implies that our fixed part contains at most
$Nl+l(l-3)/2$ points of $E$ (the same argument as in Proposition 3 c)). 
Since $l\leq 2N$, this is at most
$2Nl$. A general member of the free part, by the property (2), contains
less than $N(4N-l)$ points of $E$. So a general curve of our linear system 
contains less than
$2Nl+N(4N-l)=4N^2+Nl\leq 6N^2$ points of $E$, a contradiction.

\

{\bf Remarks:} 1) For small values of $N$, one obtains a better bound
by the same type of arguments. 
An interested reader can verify that for $N=4$, by taking sextics through
24 points of $E$ and applying Proposition 3, we obtain $|E|<24$ (and of course
for $N=3$, by taking quartics through 12 points of $E$, we obtain
$|E|<12$).

2) ({\it one of possible proofs of Theorem 2:})
 With a more elaborate argument of this kind, involving a Cayley-Bacharach
type statement, we arrive to the bound
of Theorem 2, that is, $4N^2$. This can be done as follows: suppose that
$E$ has $4N^2$ or more points, fix a subset $V$ of $4N^2$ points and
 consider curves of degree $4N$ containing $V$. This is a linear system
of dimension at least $4N^2+6N$. In the same way as above, we see that
the degree $l$ of its fixed part is less than $\frac{4}{3}N$. The fixed part
contains less than
$Nl+\frac{(l-1)(l-2)}{2}$ points of $V$,
 and so on the free part there is a set $V'$ of
at least $4N^2-Nl-\frac{l(l-3)}{2}$ of them. This is not enough to get
a contradiction immediately, but let us remark that the conditions imposed
by the points of $V'$ on curves of degree $4N-l$ are far from 
being independent. Indeed, the number of those conditions 
 is at most $\frac{(4N-l)(4N-l+3)}{2}-4N^2-6N
=4N^2-4Nl+\frac{l(l-3)}{2},$ whereas the number of points in $V'$
is at least $4N^2-Nl-\frac{l(l-3)}{2}$. When $l<\frac{4}{3}N$, the number of 
conditions is clearly inferior to the number of points: we obtain that no
collection of
$4N^2-4Nl+\frac{(l-1)(l-2)}{2}$ points
of $V'$ imposes independent conditions on curves of degree $4N-l$. 
At this point, let us recall
the following result from \cite{EP} (Corollaire 2):

{\it Let $M$ be a set of $A$ points in the plane, 
which do not impose independent 
conditions on curves of degree $\tau$. Let $s$ be a positive integer 
such that $s^2\leq A$
and $\tau>s-3+\frac{A}{s}$. Then there exists an integer $t$, $0<t<s$,
and a subset $M'\subset M$ consisting of at least $t(\tau-t+3)$ points,
which is contained in a curve of degree $t$.}

Straightforward calculations show that this result applies
to $A=4N^2-4Nl+\frac{(l-1)(l-2)}{2}$, $\tau=4N-l$ and $s=[\sqrt{A}]$, and,
moreover, that
for any $0<t<s$, we have $t(\tau-t+3)\geq Nt+\frac{(t-1)(t-2)}{2}.$
Therefore
there is a curve of degree $t$ containing at least $Nt+\frac{(t-1)(t-2)}{2}$
points of $V$. But this contradicts Proposition 3, so we are done:
$E$ contains less than $4N^2$ elements.

It is also possible to give an elementary and self-contained 
proof of Theorem 2. This is what we shall do in the next
section.

\

{\bf 5. Points of high multiplicity on plane curves}

\

Recall that, for example, in the proof of Theorem 1, it was an
essential
observation (resulting from Proposition 1) that a point of complete
ramification is a point of multiplicity at least $m^2-1$ on the
divisor
$f_*R$ of degree $3m^2-3m$. Remark that if our $R$ was irreducible, or
even reduced, this
would imply without any further work that $f$ is completely ramified
at 8 points at most. Indeed, a simple calculation shows that a reduced
curve of degree $3m^2-3m$ cannot have more than 8 points of
multiplicity $m^2-1$ or more (provided that $m>1$, of course). In
general,
let $s>1$ be an integer, and denote by $E_s(B)$ 
the set of points of multiplicity at least $s$ on a reduced plane
curve $B$ of degree $b$. The number of elements in $E_s(B)$ 
can be bounded by a function depending only
on the ratio $\frac{b}{s}$.
Indeed, a point of multiplicity at least $s$ lowers the geometric
genus
of a curve by at least $s(s-1)/2$, and the lowest possible geometric
genus of a reduced curve of degree $b$ is $1-b$, so that we get
an inequality $|E_s(B)| s(s-1)/2\leq (b-1)(b-2)/2+1-b$, or
$|E_s(B)|\leq \frac{b(b-1)}{s(s-1)},$ and we can roughly bound this, say,
by $2(\frac{b}{s})^2$.

Of course $f_*R$ is not necessarily reduced, and for an
arbitrary effective divisor
$B=\sum a_iB_i$ of degree $b$ in the plane, the observation above 
fails completely. 
One trivial reason for this is that $B$ can have a component of multiplicity
$s$,
and then $E_s(B)$ is an infinite set. But even if we suppose that it does
not, there are still some obvious counterexamples. For instance,
$B$ might be the union of a line of multiplicity $s-1$ and $b-s+1$
other
lines. Such a configuration can have $b-s+1$ points of multiplicity
$s$,
and this is not bounded by a function of $\frac{b}{s}$.

However one can ask if such a bound holds if in addition we suppose
that the low-degree components of $B$ contain ``not too many'' points
of $E_s(B)$. The following proposition provides some bound
under a very strong restriction of this kind. Together with
Proposition
3, it implies Theorem 2.

\

{\bf Proposition 5:} {\it Let $B=\sum a_iB_i$ be an effective divisor
  of
degree b
in the plane, and fix a subset $E_s$ of the set of points of multiplicity
  at least
$s$ on $B$. Set $N=[\frac{b}{s}]+1$. Suppose that 
each irreducible component $B_i$ of degree $\leq\frac{b}{s}$
contains at most $Ndeg(B_i)-1$ points of $E_s$ in its smooth
locus.  
Then $|E_s|<4N^2$.}

\

{\bf Remarks:} 1) The notation ``$B$'' comes from the ``branch divisor''.
If $B$ is $f_*R$ from the last section and $s$ is the 
lower bound for the multiplicity at $P_i$ coming from Proposition 1,
then $[\frac{b}{s}]$ is the dimension of the ambient
projective space, that is, our $N$ coincides with that of the last
section.

2) We have to consider an arbitrary subset 
of the set of points of multiplicity
  at least
$s$ on $B$, rather than this whole set, because we want to deduce
Theorem 2 from this proposition by taking the set of $P_i$'s for
$E_s$. 
Of course on
$B=f_*R$ there can be other points of high multiplicity than the
points $P_i$. The assumption of the proposition implies immediately
that $E_s$ is finite: indeed, on any component which might be of
multiplicity $s$ or more, there is only a finite number of elements
of $E_s$.

\

The proof of this proposition is based on Lemma 5 below.
Let us first introduce some more notations. Consider an irreducible
plane curve $D$ of degree $d$, and let $h$ be a function on $D\cap
E_s$, taking values in the interval $[0,1]$. For an integer $r\geq
0$, 
define $$q_{r,h}(D)=\sum_{y\in D\cap E_s} h(y)(mult_y(D))^r.$$ 

Write $B=aD+D'$, where $a\geq 0$ and $D'$ is an effective divisor not
containing $D$.

\

{\bf Lemma 5:} {\it For any $h$ as above, we have an inequality

$$bd-sq_{1,h}(D)\geq a(d^2-q_{2,h}(D)).$$}

\

{\it Proof:} The intersection number $D'D$ is equal to $(b-ad)d$.
On the other hand, 

$$ D'D\geq \sum_{y\in D\cap E_s}mult_y(D')mult_y(D)
\geq \sum_{y\in D\cap E_s}h(y)(mult_y(B)-mult_y(D)a)mult_y(D)\geq $$

$$\geq \sum_{y\in D\cap E_s}h(y)(s-mult_y(D)a)mult_y(D)=
sq_{1,h}(D)-aq_{2,h}(D),$$ and our inequality follows.

\

Now let $q_{reg}(D)$ (resp. $q_{sing}(D)$) 
denote the number of smooth (resp. singular) points of $D$ which are in $E_s$.
Let $\mu(D)$ denote the sum $\sum_{x\in Sing(D)\cap E_s}mult_x(D).$
One checks easily that $\mu(D)<d^2$ and that $q_{sing}(D)\leq
\mu(D)/2<d^2/2.$

Lemma 5 has the following long 

\

{\bf Corollary:} {\it a) If $D$ is not a component of $B$, then 
$|D\cap E_s|\leq
\frac{bd}{s}$;

b) If $D$ is a component of $B$ and $d>\frac{b}{s}$, then
$q_{reg}(D)\leq
\frac{bd}{s}$; 

c) Let $I$ be any subset of $D\cap E_s$, such that
$\sum_{x\in I}mult_x(D)\leq d^2$. Then $|I|\leq \frac{bd}{s}$.
In particular, $q_{sing}(D)\leq \frac{bd}{s}.$

d) If $q_{reg}(D)<d^2$, then 
$\frac{q_{reg}(D)}{2}+q_{sing}(D)\leq \frac{bd}{s}.$}

\

{\it Proof:} a) In this case $a=0$. The assertion follows from Lemma 5
if we take the function $h\equiv 1$ and remark that 
$q_{1,h}(D)\geq |D\cap E_s|$. 

b) Let $h$ be equal to one on smooth points of $D$ and to zero on
singular
points of $D$. Then $q_{1,h}(D)=q_{2,h}(D)=q_{reg}(D)$, and the
inequality of Lemma 5 
reads as follows: $(b-ad)d\geq (s-a)q_{reg}(D)$. Remark that by
assumption
$d>\frac{b}{s}$, so
$D$ is of multiplicity less than $s$ in $B$,
so $s-a>0$. If we suppose that 
$q_{reg}(D)>\frac{bd}{s},$ we obtain that $(b-ad)>(s-a)\frac{b}{s}$, or, since
$a\neq 0$, that $d<\frac{b}{s}$, contradicting the assumption.

c) Take $h=\frac{1}{mult_y(D)}$ if $y\in I$, and zero otherwise.
Lemma 5 gives then $bd-s|I|\geq 0$, q.e.d.

d) If $q_{reg}(D)+\mu(D)\leq d^2$, then by c), $|D\cap E_s|\leq
\frac{bd}{s}$
and the assertion is obvious, so we may assume that 
$q_{reg}(D)+\mu(D)> d^2.$ Consider the number 
$t=\frac{d^2-q_{reg}(D)}{\mu(D)}$, by assumption we have then $0<t<1$.
So Lemma 5 applies to the following function $h$:
$h(y)=\frac{t}{mult_y(D)}$
if $y$ is singular on $D$ and $h(y)=1$ if $y$ is  smooth
on $D$. The inequality of Lemma 5 rewrites as

$bd-s(q_{reg}(D)+tq_{sing}(D))\geq a(d^2-q_{reg}(D)-t\mu(D))=0,$

which gives
$$\frac{bd}{s}\geq
q_{reg}(D)+\frac{d^2-q_{reg}(D)}{\mu(D)}q_{sing}(D)=
q_{reg}(D)(1-\frac{q_{sing}(D)}{\mu(D)})+q_{sing}(D)\frac{d^2}{\mu(D)}.$$
Our assertion follows now from the remarks made just before the
Corollary:
$\mu(D)<d^2$ and $q_{sing}(D)\leq \mu(D)/2$.

\

{\it Proof of Proposition 5:} First let us show that the assumption,
together with the Corollary, implies that 
{\it any} curve $D$ of degree $d$ contains less than
$\frac{3}{2}Nd$
points of $E_s$. Indeed, we may suppose that $D$ is
irreducible; if $d>\frac{b}{s}$, then by a) and b) of the Corollary,
we have $q_{reg}(D)\leq \frac{bd}{s}<d^2$, that is, we can apply d).
And this implies that 
$|D\cap E_s|=\frac{1}{2}q_{reg}(D)+(\frac{1}{2}q_{reg}(D)
+q_{sing}(D))\leq \frac{3}{2}\frac{bd}{s}<\frac{3}{2}Nd$ as required.

If $d\leq\frac{b}{s}$, then $q_{reg}(D)<Nd$, by
assumption
if $D$ is a component of $B$ and by a), if not. It remains to remark
that $q_{sing}(D)<d^2/2\leq \frac{bd}{2s}$, so again
$|D\cap E_s|<\frac{3}{2}Nd.$

Also, recall that the Corollary says that a curve of degree $d$ 
which does not have common components with $B$, contains at most
$\frac{bd}{s}$ (so less than $Nd$) points of $E_s$. 

From these two facts we deduce $|E_s|<4N^2$: indeed, suppose that $E_s$ has
$4N^2$ or more elements. Consider the curves of degree $3N$ passing
through $4N^2$ points of $E_s$. It is a linear system of dimension
at least $\frac{3N(3N+3)}{2}-4N^2$. Let $l$ be the degree of the union
of
its base components. The degree of the free part is then $3N-l$, and
the dimension of the linear system is at most
$\frac{(3N-l)(3N-l+3)}{2}$,
so we get an inequality $N^2+9N\leq (3N-l)(3N-l+3)$, from which one
deduces $l\leq 2N$. Now the free part contains at most
$\frac{b}{s}(3N-l)$
points of $E_s$, because its general member does not contain any
component
of $B$. The fixed part contains less than $\frac{3}{2}Nl$
points of $E_s$. So we conclude that a general curve of our linear
system
contains less than $N(3N-l)+ \frac{3}{2}Nl=3N^2+\frac{1}{2}Nl\leq
4N^2$
points of $E_s$, and this is a contradiction. Proposition 5 is proved.

\

{\it Proof of Theorem 2:} Let us recall the situation:
 we have a finite
morphism $f_T: T\rightarrow \p^2$,
where $T$ is a smooth complete intersection surface
of type $(m,m,\dots,m)$ in $\p^{N-1}$, and $f^*H_{\p^2}=mH_T$ where $H$
stands for the hyperplane section divisor. We suppose that $m>>N$, so that
Proposition 3 applies. By adjunction and Hurwitz'
formula, we compute the ramification divisor $R$ of $f_T$: $R$ is
linearly equivalent to $(m-1)NH_T$. By projection formula, $f_*R$
is a divisor of degree $m^{N-2}(m-1)N$. Furthermore, we have a set
$E$ of points $P_i$, $1\leq i\leq k$, on $\p^2$, and we know
that the inverse image of each $P_i$ consists of $m^{N-3}$ points on $T$,
and in each of those, the local degree of the map $f_T$ is $m^2$.

By Proposition 1, $f_*R$ is of multiplicity at least $m^{N-3}(m^2-1)$
at $P_i$, for any $1\leq i\leq k$. It remains to apply Proposition 5
to $B=f_*R$, $b=m^{N-2}(m-1)N$, $s=m^{N-3}(m^2-1)$ and
$E_s=E$. Indeed,
Proposition 3 implies that the condition of Proposition 5 is satisfied.

{}

\end{document}